\documentclass[10pt,reqno]{amsart}
\usepackage{float}
\usepackage{graphicx}
\usepackage{graphics}
\usepackage{subfigure}
\usepackage{mathrsfs}
\usepackage[colorlinks,linkcolor=black,anchorcolor=black,citecolor=black,hyperindex=black,CJKbookmarks=black]{hyperref}

\newtheorem{theorem}{Theorem}[section]
\newtheorem*{lemma}{Lemma A}

\newtheorem{proposition}[theorem]{Proposition}

\theoremstyle{definition}
\newtheorem{definition}[theorem]{Definition}

\numberwithin{equation}{section}

\begin{document}

\title[]{A useful lemma for calculating the Hausdorff dimension of certain sets in Engel expansions}

\author {Lei Shang}

\subjclass[2010]{Primary 11K55; Secondary 28A80}
\keywords{Engel expansions, Growth rate of digits, Hausdorff dimension}

\begin{abstract}
Let $\{s_n\}$ and $\{t_n\}$ be two sequences of positive real numbers. Under some mild conditions on $\{s_n\}$ and $\{t_n\}$, we give the precise formula of the Hausdorff dimension of the set
\[
\mathbb{E}(\{s_n\},\{t_n\}):=\Big\{x\in(0,1): s_{n}<d_{n}(x)\leq s_n+t_n, \forall n\geq1\Big\},
\]
where $d_n(x)$ denotes the digit of the Engel expansion of $x$. This result improves the Lemma 2.6 of Shang and Wu (2021JNT), and is very useful for calculating the Hausdorff dimension of certain sets in Engel expansions.
\end{abstract}

\maketitle

\section{Introduction}

Let $T:[0,1)\to [0,1)$ be the \emph{Engel expansion map} defined by $T(0):=0$ and
\[
T(x) := x\left\lceil\frac{1}{x}\right\rceil -1,\ \ \ \forall x \in (0,1),
\]
where $\lceil y\rceil$ denotes the least integer not less than $y$. Denote by $T^k$ the $k$th iteration of $T$. For $x\in (0,1)$, if $x$ is rational, then there exists $n \in \mathbb{N}$ such that $T^n(x)=0$; if $x$ is irrational, then $T^n(x)>0$ for all $n \geq 1$. Then every irrational number $x\in (0,1)$ admits an infinite series with the form
\begin{equation}\label{EE}
x= \frac{1}{d_1(x)} + \frac{1}{d_1(x)d_2(x)}+\cdots+\frac{1}{d_1(x)\cdots d_n(x)}+\cdots
\end{equation}
by letting $d_1(x) := \lceil 1/x\rceil$ and $d_{n+1}(x) := d_1(T^n(x))$ for all $n \geq 1$. The expression \eqref{EE} is called the \emph{Engel expansion} of $x$ and the integers $d_n(x)$ are called the \emph{digits} of the Engel expansion of $x$. It was remarked in \cite[p.\,7]{ERS58} that $d_{n+1}(x) \geq d_n(x) \geq 2$ for all $n \geq 1$ and $d_n(x) \to \infty$ as $n \to \infty$. 
We refer the reader to Galambos \cite{Gal76} for more information of Engel expansions.

Let $\{s_n\}$ and $\{t_n\}$ be sequences of positive real numbers. Write
\[
\mathbb{E}(\{s_n\},\{t_n\}):=\Big\{x\in(0,1): s_{n}<d_{n}(x)\leq s_n+t_n, \forall n\geq1\Big\}.
\]
We will prove the following result.

\begin{lemma}\label{A}
Assume that: (1)\,$s_n \geq t_n \geq 2$; (2)\,$s_{n+1} \geq s_n+t_n$; (3)\,$\lim_{n \to \infty}s_n=\infty$. Then
\begin{align}\label{Formula}
\dim_\mathrm{H}\mathbb{E}(\{s_n\},\{t_n\})=\liminf_{n \to\infty}\frac{\sum_{k=1}^n \log t_k}{\sum_{k=1}^{n+1} \log s_k+\log s_{n+1}-\log t_{n+1}}.
\end{align}
\end{lemma}

We remark that a similar result was obtained by Liao and Rams \cite{LR} for a class of infinite iterated function systems, while the Engel expansion system does not belong to their setting, see \cite{JR, LR}.
We also point out that the result of Lemma A would be very useful for calculating the Hausdorff dimension of certain sets in Engel expansions, see for example \cite{FWnon, LW03, LL, SWjmaa}.

The paper is organized as follows. Section 2 is devoted to several definitions and basic properties of the Engel expansion. The proof of Lemma A will be given in Section 3.

\section{Preliminaries}

\begin{definition}
A finite sequence $(\sigma_1, \cdots, \sigma_n) \in \mathbb{N}^n$ is said to be \emph{admissible} if there exists $x \in (0,1)$ such that $d_k(x) = \sigma_k$ for all $1\leq k\leq n$. An infinite sequence $(\sigma_1, \cdots, \sigma_k, \cdots) \in \mathbb{N}^\mathbb{N}$ is said to be \emph{admissible} if there exists $x \in (0,1)$ such that for all $n \geq 1$, $d_k(x) = \sigma_k, \forall 1 \leq k \leq n$.

\end{definition}

Denote by $\Sigma_n$ the collection of all admissible sequences with length $n$ and by $\Sigma$ that of all infinite admissible sequences.
The following result gives a characterisation of admissible sequences.

\begin{proposition}[\cite{Gal76}]\label{AD}
$(\sigma_1, \cdots, \sigma_n) \in \Sigma_n$ if and only if $2 \leq \sigma_1  \leq \cdots \leq \sigma_n$;
$(\sigma_1,  \cdots, \sigma_n, \cdots) \in \Sigma$ if and only if
\[
\sigma_{n+1} \geq \sigma_n \geq 2, \ \forall n \geq 1 \ \ \ \ \ \text{and} \ \ \ \ \ \lim_{n \to \infty}\sigma_n = \infty.
\]
\end{proposition}

\begin{definition}
Let $(\sigma_1, \cdots, \sigma_n) \in \Sigma_n$. We call
\begin{equation*}
I_n(\sigma_1, \cdots, \sigma_n) := \big\{x \in (0,1): d_1(x)=\sigma_1, \cdots,d_n(x)=\sigma_n\big\}
\end{equation*}
the \emph{cylinder} of order $n$ associated to $(\sigma_1,\dots,\sigma_n)$.
\end{definition}

We use $|I|$ to denote the diameter of an interval $I$.

\begin{proposition}[{\cite{Gal76}}]\label{cylinder}
Let $(\sigma_1, \cdots, \sigma_n) \in \Sigma_n$. Then $I_n(\sigma_1, \cdots, \sigma_n) =[A_n, B_n)$, where
\[
A_n:=\frac{1}{\sigma_1}+\cdots+\frac{1}{\sigma_1\sigma_2\cdots \sigma_{n-1}}+ \frac{1}{\sigma_1\sigma_2\cdots \sigma_{n-1}\sigma_n}
\]
and
\[
B:=\frac{1}{\sigma_1}+\cdots+ \frac{1}{\sigma_1\sigma_2\cdots \sigma_{n-1}}+ \frac{1}{\sigma_1\sigma_2\cdots \sigma_{n-1}(\sigma_n-1)}.
\]
Moreover,
\begin{equation*}\label{cylinder length}
\left|I_n(\sigma_1, \sigma_2, \cdots, \sigma_n)\right| = \frac{1}{\sigma_1\sigma_2\cdots \sigma_{n-1}\sigma_n(\sigma_n-1)}.
\end{equation*}
\end{proposition}

\begin{proposition}[{\cite[Proposition 4.1]{Fal90}}]\label{upp}
Suppose $\mathbb{F}$ can be covered by $\mathcal{N}_n$ sets of diameter at most $\delta_n$ with $\delta_n \to 0$ as $n \to \infty$. Then
\[
\dim_{\rm H}\mathbb{F} \leq \liminf_{n \to \infty} \frac{\log \mathcal{N}_n}{-\log \delta_n}.
\]
\end{proposition}

\begin{proposition}[{\cite[Example 4.6]{Fal90}}]\label{low}
Let $[0,1] = \mathbb{E}_0 \supset \mathbb{E}_1 \supset \cdots$ be a decreasing sequence of sets and $\mathbb{E} = \bigcap_{n \geq 0} \mathbb{E}_n$. Assume that each $\mathbb{E}_n$ is a union of a finite number of disjoint closed intervals (called basic intervals of order $n$) and each basic interval in $\mathbb{E}_{n-1}$ contains $m_n$ intervals of $\mathbb{E}_n$ which are separated by gaps of lengths at least $\varepsilon_n$.
If $m_n \geq 2$ and $\varepsilon_{n-1}> \varepsilon_n >0$, then
\[
\dim_\mathrm{H} \mathbb{E} \geq \liminf_{n \to \infty} \frac{\log(m_1m_2 \cdots m_{n-1})}{-\log (m_{n}\varepsilon_n)}.
\]
\end{proposition}

\section{Proof of Lemma A}

In this section, we will give the proof of Lemma A. Assume that: (1)\,$s_n \geq t_n \geq 2$; (2)\,$s_{n+1} \geq s_n+t_n$; (3)\,$\lim_{n \to \infty}s_n=\infty$. Let
\begin{equation*}
\mathcal{D}_n:=\big\{(\sigma_1,\cdots \sigma_n)\in \mathbb{N}^n:s_k<\sigma_k\leq s_k+t_k,\forall 1\leq k\leq n\big\}.
\end{equation*}
Note that $\sigma_1>s_1\geq 2$ and $\sigma_{k+1}>s_{k+1} \geq s_k+t_k >\sigma_k$, so $\sigma_n \geq \cdots\geq \sigma_1\geq 2$. That is to say, $(\sigma_1,\cdots \sigma_n)$ is admissible. For $(\sigma_1,\cdots \sigma_n) \in \mathcal{D}_n$, let
\begin{equation*}
J_n(\sigma_1,\cdots,\sigma_n):=\bigcup_{s_{n+1}<j\leq s_{n+1}+t_{n+1}} cl(I_{n+1}(\sigma_1,\cdots,\sigma_n,j)),
\end{equation*}
where $cl(\cdot)$ denotes the closure of a set. By Proposition \ref{cylinder}, we know that $J_n(\sigma_1,\cdots,\sigma_n)$ is a closed interval, which is called the \emph{basic interval} of order $n$. Write
\begin{equation}\label{En}
\mathbb{E}_n:=\bigcup_{(\sigma_1,\cdots,\sigma_n)\in \mathcal{D}_n} J_n(\sigma_1,\cdots,\sigma_n)
\end{equation}
with the convention $\mathbb{E}_0:=[0,1]$. Then
\begin{equation*}
\mathbb{E}(\{s_n\},\{t_n\})=\bigcap_{n=0}^\infty\mathbb{E}_n.
\end{equation*}

For the upper bound of $\dim_{\rm H}\mathbb{E}(\{s_n\},\{t_n\})$, for any $n \geq 1$, we derive from \eqref{En} that $\{J_n(\sigma_1,\cdots,\sigma_n):(\sigma_1,\cdots,\sigma_n) \in \mathcal{D}_n\}$ is a cover of $\mathbb{E}(\{s_n\},\{t_n\})$. Then
\begin{align*}
\mathcal{N}_n:=\#\mathcal{D}_n=(\lfloor s_1+t_1\rfloor-\lfloor s_1\rfloor)\cdots(\lfloor s_n+t_n\rfloor-\lfloor s_n\rfloor)\leq 2^nt_1\cdots t_n.
\end{align*}
For any $(\sigma_1,\cdots \sigma_n) \in  \mathcal{D}_n$, by Proposition \ref{cylinder}, we see that
\begin{align*}
|J_{n}(\sigma_1,\cdots,\sigma_n)|&=\sum_{j=\lfloor s_{n+1}\rfloor+1}^{\lfloor s_{n+1}+t_{n+1}\rfloor}|I(\sigma_1,\cdots,\sigma_n,j)|
 =\frac{1}{\sigma_1\cdots \sigma_n}\sum_{\lfloor s_{n+1}\rfloor+1}^{\lfloor s_{n+1}+t_{n+1}\rfloor}\frac{1}{j(j-1)}.
\end{align*}
Note that $s_{n+1}-1\geq s_{n+1}/2$ and $s_{n+1}+t_{n+1}\geq s_{n+1}$, so
\begin{align*}
\sum_{j=\lfloor s_{n+1}\rfloor+1}^{\lfloor s_{n+1}+t_{n+1}\rfloor}\frac{1}{j(j-1)} = \frac{1}{\lfloor s_{n+1}\rfloor} - \frac{1}{\lfloor s_{n+1}+t_{n+1}\rfloor}
 \leq \frac{1}{s_{n+1}-1}-\frac{1}{s_{n+1}+t_{n+1}}
\leq\frac{4t_{n+1}}{s_{n+1}^2}.
\end{align*}
Since $\sigma_k\geq s_k$, we have
\begin{align*}
|J_n(\sigma_1,\cdots,\sigma_n)|\leq\frac{1}{s_1\cdots s_n}\cdot\frac{4t_{n+1}}{s_{n+1}^2}=:\delta_n.
\end{align*}
From Proposition \ref{upp}, we conclude that
\begin{align*}
\dim_\mathrm{H}\mathbb{E}(\{s_n\},\{t_n\})&\leq \liminf_{n\to \infty}\frac{\log\mathcal{N}_n}{-\log \delta_n}\\
&\leq\liminf_{n\to \infty}\frac{n\log2+\log(t_1\cdots t_n)}{\log(s_1\cdots s_{n+1})+\log s_{n+1}-\log t_{n+1}-\log 4}\\
&=\liminf_{n\to \infty}\frac{\sum_{k=1}^n\log t_{k}}{\sum_{k=1}^{n+1}\log s_{k}+\log s_{n+1}-\log t_{n+1}}.
\end{align*}

For the lower bound of $\dim_{\rm H}\mathbb{E}(\{s_n\},\{t_n\})$, by the structure of basic intervals, we deduce that each basic interval of order $n-1$ contains
\[
\frac{t_n}{2}<\lfloor t_n\rfloor\leq m_n:=\lfloor s_n+t_n\rfloor - \lfloor s_n\rfloor < t_n+1 <2t_n
\]
basic intervals of order $n$. Next we will estimate the gaps between two basic intervals with the same order. For two adjacent sequences $(\sigma_1,\cdots,\sigma_n)$ and $(\sigma^{\prime}_1,\cdots,\sigma^{\prime}_n)$ in $\mathcal{D}_n$, we deduce that the cylinders $J_n(\sigma_1,\cdots,\sigma_n)$ and $J_n(\sigma^{\prime}_1,\cdots,\sigma^{\prime}_n)$ are different and then $J_n(\sigma_1,\cdots,\sigma_n)$ is on the left-hand or right-hand side of $J_n(\sigma^{\prime}_1,\cdots,\sigma^{\prime}_n)$. Without loss of generality, we assume that $J_n(\sigma_1,\cdots,\sigma_n)$ is on the left-hand side of $J_n(\sigma^{\prime}_1,\cdots,\sigma^{\prime}_n)$. Then
\[
\Pi_1:=\bigcup_{\sigma_n\leq j\leq \lfloor s_{n+1}\rfloor}I_{n+1}(\sigma_1,\cdots,\sigma_n,j) \ \text{and}\ \Pi_2:=\bigcup_{j\geq \lfloor s_{n+1}+t_{n+1}\rfloor+1}I_{n+1}(\sigma^{\prime}_1,\cdots,\sigma^{\prime}_n,j)
\]
are the gap of $J_n(\sigma_1,\cdots,\sigma_n)$ and $J_n(\sigma^{\prime}_1,\cdots,\sigma^{\prime}_n)$. For $\Pi_1$, we have
\begin{align*}
|\Pi_1|&=\sum_{\sigma_n\leq j\leq \lfloor s_{n+1}\rfloor}|I_{n+1}(\sigma_1,\cdots,\sigma_n,j)|\\
&=\frac{1}{\sigma_1\cdots\sigma_n}\sum_{\sigma_n\leq j\leq \lfloor s_{n+1}\rfloor}\frac{1}{j(j-1)}\\
&=\frac{1}{\sigma_1\cdots\sigma_n}\left(\frac{1}{\sigma_n-1}-\frac{1}{\lfloor s_{n+1}\rfloor}\right)\\
&\geq\frac{1}{2s_1\cdots 2s_n}\left(\frac{1}{s_n+t_n-1}-\frac{1}{s_{n+1}-1}\right).
\end{align*}
For $\Pi_2$, we obtain
\begin{align*}
|\Pi_2|&=\sum_{j\geq \lfloor s_{n+1}+t_{n+1}\rfloor+1}|I_{n+1}(\sigma^{\prime}_1,\cdots,\sigma^{\prime}_n,j)|\\
&=\frac{1}{\sigma^{\prime}_1\cdots\sigma^{\prime}_n}\sum_{j\geq \lfloor s_{n+1}+t_{n+1}\rfloor+1}\frac{1}{j(j-1)}\\
&\geq\frac{1}{2s_1\cdots 2s_n}\cdot\frac{1}{s_{n+1}+t_{n+1}}.
\end{align*}
Note that $s_{n+1} \geq s_n+t_n$, so
\begin{align*}
\frac{1}{s_n+t_n-1}-\frac{1}{s_{n+1}-1} +\frac{1}{s_{n+1}+t_{n+1}} &= \frac{1}{s_n+t_n-1}- \frac{t_{n+1}+1}{(s_{n+1}-1)(s_{n+1}+t_{n+1})}\\
&\geq \frac{1}{s_n+t_n-1}\cdot \left(1- \frac{t_{n+1}+1}{s_{n+1}+t_{n+1}}\right)\\
&=\frac{1}{s_n+t_n-1}\cdot\frac{s_{n+1}-1}{s_{n+1}+t_{n+1}}.
\end{align*}
Since $s_n+t_n-1 < 2s_n$, $s_{n+1}-1 \geq s_{n+1}/2$ and $s_{n+1}+t_{n+1} \leq 2s_{n+1}$, we see that the length of the gap of $J(\sigma_1,\cdots,\sigma_n)$ and $J(\sigma^{\prime}_1,\cdots,\sigma^{\prime}_n)$ is at least
\begin{align*}
|\Pi_1|+|\Pi_2|
&\geq\frac{1}{2s_1\cdots2s_n}\cdot\frac{1}{s_n+t_n-1}\cdot\frac{s_{n+1}-1}{s_{n+1}+t_{n+1}}\\
&\geq \frac{1}{2s_1\cdots2s_n}\cdot\frac{1}{8s_n}\\
&=\frac{1}{2^{n+3}}\cdot\frac{1}{s_1\cdots s_ns_n}:=\varepsilon_n.
\end{align*}
It follows from Proposition \ref{low} that
\begin{align*}
\dim_\mathrm{H}\mathbb{E}(\{s_n\},\{t_n\})&\geq \liminf_{n\to \infty}\frac{\log(m_1m_2\cdots m_{n-1})}{-\log(m_n\varepsilon_n)}\\
&\geq\liminf_{n\to \infty}\frac{-n\log2+\sum _{k=1}^n \log t_k }{(n+5)\log2+\sum_{k=1}^{n+1}\log s_k+\log s_{n+1}-\log t_{n+1}}\\
&=\liminf_{n\to \infty}\frac{\sum _{k=1}^n \log t_k}{\sum_{k=1}^{n+1}\log s_k+\log s_{n+1}-\log t_{n+1}}.
\end{align*}


\begin{thebibliography}{10}
\bibitem{ERS58} P. Erd\H{o}s, A. R\'{e}nyi and P. Sz\"{u}sz, {\it On Engel's and Sylvester's series}, Ann. Univ. Sci. Budapest. E\"{o}tv\"{o}s. Sect. Math. 1 (1958), 7--32.

\bibitem{Fal90} K. Falconer, {\it Fractal Geometry: Mathematical Foundations and Applications}, John Wiley \& Sons, Ltd., Chichester, 1990.



\bibitem{FWnon} L. Fang and M. Wu, {\it Hausdorff dimension of certain sets arising in Engel expansions}, Nonlinearity 31 (2018), 2105--2125.




\bibitem{Gal76} J. Galambos, {\it Representations of Real Numbers by Infinite Series}, Lecture Notes in Mathematics, Vol. 502. Springer-Verlag, Berlin-New York, 1976.


\bibitem{JR} T. Jordan and M. Rams, {\it Increasing digit subsystems of infinite iterated function systems}, Proc. Amer. Math. Soc. 140 (2012), 1267--1279.


\bibitem{LR} L. Liao and M. Rams, {\it Big Birkhoff sums in $d$-decaying Gauss like iterated function systems}, to appear in Studia Mathematica.




\bibitem{LW03} Y. Liu and J. Wu, {\it Some exceptional sets in Engel expansions}, Nonlinearity 16 (2003), 559--566.


\bibitem{LL} M. L\"{u} and J. Liu, {\it Hausdorff dimensions of some exceptional sets in Engel expansions}, J. Number Theory 185 (2018), 490--498.





\bibitem{SWjnt} L. Shang and M. Wu, {\it On the growth speed of digits in Engel expansions}, J. Number Theory 219 (2021), 368--385.


\bibitem{SWjmaa} L. Shang and M. Wu, {\it On the exponent of convergence of the digit sequence of Engel series}, J. Math. Anal. Appl. 504 (2021), Paper No. 125368.
















\end{thebibliography}
\end{document}